\documentclass[10pt,draftcls,onecolumn]{IEEEtran}
\usepackage{epsfig,amssymb,latexsym,color,pifont,colordvi,multicol}
\usepackage[cmex10]{amsmath}
\usepackage{subfigure}

\definecolor{backgrey}{rgb}{0.86,0.86,0.86}
\definecolor{dblue}{rgb}{0,0.0,0.5}
\definecolor{dred}{rgb}{0.4,0.2,0}
\definecolor{dgreen}{rgb}{0.0,0.5,0}

\newcommand{\captionfonts}{\small}
\makeatletter  
\long\def\@makecaption#1#2{%
  \vskip\abovecaptionskip
  \sbox\@tempboxa{{\captionfonts #1: #2}}%
  \ifdim \wd\@tempboxa >\hsize
    {\captionfonts #1: #2\par}
  \else
    \hbox to\hsize{\hfil\box\@tempboxa\hfil}%
  \fi
  \vskip\belowcaptionskip}
\makeatother   
\newtheorem{theorem}{Theorem}

\newtheorem{assumption}[theorem]{Assumption}

\newtheorem{remark}[theorem]{Remark}

\newtheorem{corollary}[theorem]{Corollary}

\newtheorem{definition}[theorem]{Definition}

\newtheorem{lemma}[theorem]{Lemma}

\title{Limitations for nonlinear observation over erasure channel}
\author{\quad Amit Diwadkar \quad Umesh Vaidya 
\thanks{A. Diwadkar is a graduate student with the Department of Electrical and
Computer Engineering, Iowa State University, Ames, IA, 50011
diwadkar@iastate.edu}
\thanks{U. Vaidya is with the Department of Electrical and
Computer Engineering,
Iowa State University,
Ames, IA, 50011 ugvaidya@iastate.edu}%
}
\begin{document}
\maketitle \thispagestyle{empty} \pagestyle{empty}

\begin{abstract}
In this paper, we study the problem of state observation of nonlinear systems over an erasure channel. The notion of mean square exponential stability is used to analyze the stability property of observer error dynamics. The main results of this paper prove, fundamental limitation arises for mean square exponential stabilization of the observer error dynamics, expressed in terms of probability of erasure, and positive Lyapunov exponents of the system. Positive Lyapunov exponents are a measure of average expansion of nearby trajectories on an attractor set for nonlinear systems. Hence, the dependence of limitation results on the Lyapunov exponents highlights the important role played by non-equilibrium dynamics in observation over an erasure channel.  The limitation on observation is also related to measure-theoretic entropy of the system, which is another measure of dynamical complexity. The limitation result for the observation of linear systems is obtained as a special case, where Lyapunov exponents are shown to  emerge as the natural generalization of eigenvalues from linear systems to nonlinear systems.
\end{abstract}


\section{Introduction}
The problem of state estimation of systems over erasure channels has attracted a lot of attention lately, given the importance of this problem in the control of systems over a network \cite{networksystems_specialissue}. The problem of state estimation with intermittent observation was first studied in \cite{nahi_estimation, estimation_hadidi}. In \cite{Luca03kalmanfiltering, murray_epstein}, state estimation over an erasure channel with  different performance metrics on the error covariance is studied. In \cite{Luca03kalmanfiltering}, under some assumptions on system dynamics, it is proved that there exists a critical non-erasure probability below which the error covariance is unbounded. A Markov jump linear system framework is used to model the state estimation problem with intermittent measurement and to provide conditions for the convergence of error covariance in \cite{estimation_MJLS}. In \cite{seiler_estimation}, state estimation over erasure channel with Markovian packet loss is studied. However, all the above results are developed for linear time invariant (LTI) systems. There is no systematic result that addresses the state estimation problem for nonlinear systems over erasure channels. Thus there is a need for extension and development of such results for nonlinear systems, with regard to their applications in network systems consisting of nonlinear components, such as power system networks, biological networks, and Internet communication networks.

In this paper, we study the problem of state observation of nonlinear systems over an erasure channel, with the objective to develop limitation results for state observation. We expect the limitation results for the state observation problem, to provide useful insight into the more challenging problem of state estimation over an erasure channel. The erasure channel is modeled as  an on/off Bernoulli switch.  We use mean square exponential (MSE) stability to study the state observation problem over an erasure channel. The main result of this paper shows, that a fundamental limitation arises in MSE stabilization of the observer error dynamics. This limitation is expressed in terms of erasure probability and {\it global} instability of the nonlinear system. In particular, under a certain ergodictiy assumption, we show  the instability of a nonlinear system can be expressed in terms of the sum of positive Lyapunov exponents of the system. Using  Ruelle's inequality from ergodic theory of a dynamical system \cite{walter_ergodic_theory}, the sum of the positive Lyapunov exponents can be related to the entropy of  a nonlinear system. Hence, the limitation result can be interpreted in terms of the entropy of a nonlinear system.  Our result involving Lyapunov exponents of a non-trivial (other than equilibrium point) invariant measure is also the first to highlight the important role played by the non-equilibrium dynamics in the limitations on nonlinear observation.

 There are two main contributions of this paper. First, it adopts and extends the formalism from erogodic theory of random dynamical systems to study the problem of nonlinear observation over an erasure channel. Second, the result  provides an analytical relationship between the maximum tolerable channel uncertainty (i.e., the maximum erasure probability) and the inability of the system to maintain mean square exponential stability of the observer error dynamics.

The organization of this paper is as follows. In section \ref{prelim}, we discuss the problem and provide necessary assumptions and stability definition. In section \ref{perf}, we prove the main results of this paper. A simulation example is presented in section \ref{sim}, followed by conclusions in section \ref{inf}.
%
\section{Preliminaries}
\label{prelim}
\noindent The set-up for nonlinear observations with a unique erasure channel at the output is described by the following equations:
\begin{align}
\label{dynamical_sys}
x_{t+1} = f(x_t),\quad \quad y_t = \xi_t h(x_t),
\end{align}
where $x_t\in X \subseteq \mathbb{R}^N$ is the state, $y_t\in Y\subseteq \mathbb{R}^M$ is the output, and $\xi_t\in \{0,1\}$ is a Bernoulli random variable with probability distribution $\text{Prob}(\xi_t=1)=p$ for all $t\geq 0$, with $0<p<1$, and independent of $\xi_{\tau}$ for $\tau\neq t$.
The IID (independent identically distributed) random variable, $\xi_t$, models the erasure channel between the plant and the observer through which all the outputs are sent to the observer simultaneously.
\begin{remark}\label{remark_unstable} To make the problem interesting, we assume that $M<N$ and $0<p<1$. The $0<p$
 assumption implies that the system dynamics, $x_{t+1}=f(x_t)$, is unstable and hence requires some non-zero probability of erasure for the observer to work.
\end{remark}
We now provide the following definition of an observability rank condition for nonlinear systems \cite{nijmeijer_observer}.
\begin{definition}[Observability Rank Condition]\label{observability_condition}
Consider the map $\theta^{N-1}(x):X\to \underbrace {Y\times\ldots \times Y}_N$
\begin{align}
\label{observability_map}
\theta^{N-1}(x) := \left(h(x),h(f(x)),\ldots,h(f^{N-1}(x)\right)'.
\end{align}
The system (\ref{dynamical_sys}) is said to satisfy the observability rank condition at $x$, if
 \[rank \left(\frac{\partial \theta^{N-1}(x)}{\partial x}\right)=N.\]
\end{definition}
We make following assumption on the system dynamics.
\begin{assumption}\label{assume}
The system mapping, $f$, and output function, $h$, are $C^r$ functions of $x$, for $r\geq 1$, with $f(0)=0$, $h(0)=0$, and the Jacobian $\frac{\partial f}{\partial x}(x)$ is uniformly bounded above and below for all $x \in X$.
Furthermore, the system satisfies the observability rank condition (Definition \ref{observability_condition}) and there exist $\alpha_{\theta} > 0$ and $\beta_{\theta}>0$, such that
\begin{align}
\label{observability_matrix}
\alpha_{\theta}I_N<\frac{\partial \theta^{N-1}}{\partial x}'(x)\frac{\partial \theta^{N-1}}{\partial x}(x) <\beta_{\theta}I_N
\end{align}
for all $x \in X$ and, $I_N$ is the $N\times N$ Identity matrix.
\end{assumption}
\begin{remark} Assumption \ref{assume} and in particular the observability rank condition are essential for the  observer design for the system with no erasure at the output.
\end{remark}

The stochastic notion of stability  we use to analyze the observer error dynamics is defined in the context of a general random dynamical system (RDS) of the form $x_{t+1}=S(x_t,\zeta_t)$,
where $x_t \in X\subseteq  \mathbb{R}^N$, $\zeta_t\in W=\{0,1\}$ for $t \geq 0$, are IID random variables with probability distribution $\text{Prob}(\zeta_t=1)=p$. The system mapping $S: X\times W\rightarrow X$ is assumed to be at least $C^1$ with respect to $x_t\in X$ and measurable w.r.t  $\zeta_t$. We assume  $x=0$ is an equilibrium point, i.e., $S(0,\zeta_t)=0$. The following notion of stability can be defined for RDS \cite{Hasminskii_book,applebaum}.
\begin{definition}[ Mean Square  Exponential (MSE) Stable]\label{mean_square_stable_def} The solution, $x=0$, is said to be MSE stable  for $x_{t+1} = S(x_t,\zeta_t)$, if there exist  positive constants $L<\infty$ and $\beta<1$, such that
\[E_{\zeta_0^t}\left [\parallel x_{t+1}\parallel^2 \right]\leq L\beta^t\|x_0\|^2,\;\; \forall t\geq 0\]
for Lebesgue almost all initial condition, $x_0\in X$, where $E_{\zeta_0^t}[\cdot]$ is the expectation taken over the sequence $\{\zeta_0,\ldots,\zeta_t\}$.
\end{definition}
\section{Main results}
\label{perf}
The main results of this paper are derived under the following assumption on the observer dynamics.
\begin{assumption}\label{observer_assumption}
The observer gain, $K$, is assumed  deterministic and not an explicit function of the channel erasure state $\xi_t$ nor its history (i.e., $\xi_0^{t-1}$). The observer dynamics is assumed to be of the form:
\begin{eqnarray}\label{observer_sys}
\hat x_{t+1} &=& f(\hat x_t) + K(y_t) - K(\hat y_t),\;\;\;\;\hat y_t=\xi_t h(\hat x_t),
\end{eqnarray}
where $\hat x\in X$ is the observer state, $\hat y\in Y$ is the observer output, and $K: Y\to X$ is the observer gain and assumed to be a $C^r$ function of $y$, for $r\geq 1$, and satisfies $K(0)=0$. Thus the property $K(0) = 0$ and $\xi_t \in \{0,1\}$, allows us to rewrite the observer dynamics (\ref{observer_sys}) as follows:
\begin{align}\label{observer_sys2}
\hat x_{t+1} = f(\hat x_t) + \xi_t K(h(x_t)) - \xi_tK(h(\hat x_t)).
\end{align}
\end{assumption}
We assume that the observer output $\hat y_t$ is an explicit function of channel state, $\xi_t$. This assumption is justified by assuming a TCP-like protocol, where the observer receives an immediate acknowledgement of the channel erasure state \cite{Luca03kalmanfiltering}.
\begin{remark}
In \cite{Luca03kalmanfiltering}, the problem of state estimation for an LTI system over an erasure channel is studied. The optimal estimator gain that minimizes the error covariance is shown to be a function of the channel erasure state history. With the estimator gain, a function of the channel erasure state history, the results in \cite{Luca03kalmanfiltering} only prove the error covariance will remain bounded and not converge to a steady state value, unlike the regular Kalman filtering problem for an LTI system with no loss of measurement. Hence, we conjecture (Assumption \ref{observer_assumption}) on the observer gain, not being a function of the channel erasure state or its history, is necessary for the error dynamics to be MSE stable.
\end{remark}
We first prove  Lemma \ref{linear_error_MSS} that provides a necessary condition for  MSE stability of the error dynamics  $x_t-\hat x_t$ in terms of MSE stability of the linearized error dynamics.
\begin{lemma}
\label{linear_error_MSS} Consider the observer dynamics in Eq. (\ref{observer_sys2}) and let the error dynamics (i.e., $e_t=x_t-\hat x_t$) be MSE stable (Definition \ref{mean_square_stable_def}).
Then, the following linearized error dynamics, $\eta_t\in \mathbb{R}^N$,
\begin{align}
\eta_{t+1}=\left(\frac{\partial f}{\partial x}(x_t)-\xi_t\frac{\partial K\circ h}{\partial x}(x_t)\right)\eta_t,\;\;\;\;\;x_{t+1}=f(x_t)
\label{linearized_error_sys}
\end{align}
is also MSE stable, i.e., there exist positive constants $L<\infty$ and $\beta<1$, such that $E_{\xi_0^t}\left[\left\| \eta_{t+1}\right\|^2\right]\leq L \beta^t\left\|\eta_0\right\|^2\;\;\;\forall t\geq 0$.
The functions $K$ and $h$ in (\ref{linearized_error_sys}) are the observer gain and output function, respectively, from Eq. (\ref{observer_sys}).
\end{lemma}
\begin{IEEEproof}
Define $g(x_t,\xi_t):= f(x_t) - \xi_tK(h(x_t))$ and ${\cal A}(x_t,\xi_t):=\frac{\partial g}{\partial x}(x_t,\xi_t)$. Then using Mean Value Theorem for the vector valued function, the error dynamics, can be written as
\begin{align}
e_{t+1} &= g(x_t,\xi_t) - g(x_t-e_t,\xi_t)=\left(\int_0^1 \frac{\partial g}{\partial x}(x_t-se_t,\xi_t) ds \right) e_t = \prod_{k=0}^t \left( \int_0^1 \mathcal{A}(x_k-se_k,\xi_k)ds \right) e_0 , \nonumber
\end{align}
Here $e_t$ is an implicit function of the initial error $e_0$, initial state $x_0$, and the sequence of uncertainties $\xi_0^{t-1}$. We define ${\cal B}_k(x_0,\xi_0^k,e_0) := \int_0^1 \mathcal{A}(x_k-se_k,\xi_k)ds$ and
${\cal B}_0^t(x_0,\xi_0^t,e_0) :=  \prod_{k=0}^t {\cal B}_k(x_0,\xi_0^k,e_0)$.
This gives
\begin{align*}
E_{\xi_0^t}\left[\parallel e_{t+1} \parallel^2\right] &= E\left[e_{t+1}'e_{t+1}\right] = e_0' E_{\xi_0^t}\left[{\cal B}_0^t(x_0,\xi_0^t,e_0)'{\cal B}_0^t(x_0,\xi_0^t,e_0)\right] e_0.
\end{align*}
Using Assumption \ref{assume}, we know there exists a positive constant  $\bar L<\infty$, such that $\parallel{\cal B}_k(x_0,\xi_0^k,\alpha e_0)\parallel< \bar L$ for Lebesgue almost all $x_0 \in X$
 and for some scalar, $\alpha>0$.
Let ${\cal B}_k(x_0,\xi_0^k,\alpha e_0)_{ij}$ denote the $i^{th}$ row $j^{th}$ column entry in ${\cal B}_k(x_0,\xi_0^k,\alpha e_0)$.
Now consider a sequence, $\{\alpha_l\}_{l = 1}^{\infty}$, such that $\lim_{l\to \infty}\alpha_l = 0$. Then, we have by Dominated Convergence Theorem \cite{Folland} and continuity of $\mathcal{A}(x_k-se_k,\xi_k)$, $\lim_{l \to \infty}{\cal B}_k(x_0,\xi_0^k,\alpha_l e_0)_{ij} = {\cal B}_k(x_0,\xi_0^k,0)_{ij}$ which implies $\lim_{l \to \infty}{\cal B}_k(x_0,\xi_0^k,\alpha_l e_0) = {\cal B}_k(x_0,\xi_0^k,0)$. Hence, we have
\begin{align}
\label{matrix_prod}
\lim_{l \to \infty}{\cal B}_0^k(x_0,\xi_0^k,\alpha_l e_0) = {\cal B}_0^k(x_0,\xi_0^k,0).
\end{align}
From MSE stability of the error, we obtain
$
 e_0' E_{\xi_0^t}\left[{\cal B}_0^t(x_0,\xi_0^t,e_0)'{\cal B}_0^t(x_0,\xi_0^t,e_0)\right] e_0 \leq  L\beta^t e_0'e_0$,
for some positive constants $L < \infty$ and $\beta<1$. Since the above inequality is true for any initial error, this will be true if the initial error vector used to compute the product of matrices is scaled by $\alpha_l$, where $\lim_{l\to \infty}\alpha_l = 0$. Substituting $\alpha_l e_0$ for $e_0$, we can write
\[e_0' E_{\xi_0^t}\left[{\cal B}_0^t(x_0,\xi_0^t,\alpha_l e_0)'{\cal B}_0^t(x_0,\xi_0^t,\alpha_l e_0)\right] e_0 \leq L\beta^t e_0'e_0.\]
Now, letting $l \to \infty$ and by Fatou's Lemma, we have
\begin{align}
\label{fatou_ineq}
e_0' E_{\xi_0^t}\left[\lim_{l \to \infty}{\cal B}_0^t(x_0,\xi_0^t,\alpha_l e_0)'{\cal B}_0^t(x_0,\xi_0^t,\alpha_l e_0)\right] e_0 &\leq \lim_{l \to \infty}e_0' E_{\xi_0^t}\left[{\cal B}_0^t(x_0,\xi_0^t,\alpha_l e_0)'{\cal B}_0^t(x_0,\xi_0^t,\alpha_l e_0)\right] e_0 \nonumber \\
&\leq L\beta^t e_0'e_0.
\end{align}
Thus, using (\ref{matrix_prod}) and (\ref{fatou_ineq}), we obtain
$e_0' E_{\xi_0^t}\left[{\cal B}_0^t(x_0,\xi_0^t,0)'{\cal B}_0^t(x_0,\xi_0^t,0)\right] e_0
\leq L\beta^t e_0'e_0$,
where ${\cal B}_0^t(x_0,\xi_0^t,0)$ is the product of the Jacobian matrices $\mathcal{A}(x_t,\xi_t)$, with zero initial error and computed along the nominal trajectory, $x_{t+1} = f(x_t)$. Hence,
\begin{align*}
E_{\xi_0^t}\left[e_0' \left(\prod_{k=0}^t\mathcal{A}(x_k,\xi_k)\right)'\left(\prod_{k=0}^t\mathcal{A}(x_k,\xi_k)\right)e_0\right]
\leq L\beta^t e_0'e_0.
\end{align*}
Since the matrices in the above equation are independent of $e_0$, we can substitute $\eta_0$ for $e_0$.
Now, using the evolution of $\eta_t$ from  Eq. (\ref{linearized_error_sys}), we obtain the desired result.
\end{IEEEproof}
Our next theorem provides the necessary condition for MSE stability of the linearized error dynamics.
\begin{theorem}
\label{lya_theorem}
 Let the $\eta_t$ dynamics for the system (\ref{linearized_error_sys}) be MSE stable (Definition \ref{mean_square_stable_def}). Then, there exists a matrix function of $x_t$, $P(x_t)$, such that $\gamma_1 I \leq P(x_t) \leq \gamma_2 I$ and
\begin{equation}
E_{\xi_t}\left[\mathcal{A}'(x_t,\xi_t)P(x_{t+1})\mathcal{A}(x_t,\xi_t)\right] < P(x_t),
\label{lyap_2mom}
\end{equation}
 for some positive constants $\gamma_1$, $\gamma_2$,  where $x_{t+1}=f(x_t)$ and ${\cal A}(x_t,\xi_t)=\frac{\partial f}{\partial x}(x_t)-\xi_t\frac{\partial K}{\partial y}(h(x_t))\frac{\partial h}{\partial x}(x_t)$ from (\ref{linearized_error_sys}).
\end{theorem}
\begin{IEEEproof}
To prove the necessary part, assume the system is MSE stable and consider the following construction of $P(x_t)$.
\begin{align*}
P(x_t) = \sum_{k=t}^{\infty}E_{\xi_t^{k}}\left[ \left(\prod_{j=t}^k\mathcal{A}(x_j,\xi_j)\right)' \left(\prod_{j=t}^k\mathcal{A}(x_j,\xi_j)\right) \right],
\end{align*}
where $E_{\xi_i^j}[\cdot] $ is the expectation over the random sequence $\{\xi_i,\ldots,\xi_j\}$. The existence of positive constants $\gamma_1, \gamma_2$ follows from the fact that $\eta_t$ dynamics is MSE stable and the Jacobian $\frac{\partial f}{\partial x}$ is bounded from above and below. The inequality (\ref{lyap_2mom}) follows from the construction of $P(x_t)$.
\end{IEEEproof}
We have Corollary \ref{corollary_MSS} to the Theorem \ref{lya_theorem}.
\begin{corollary}
\label{corollary_MSS}
Let the RDS (\ref{linearized_error_sys}) be MSE stable. Then, there exists a matrix function of $x_t$, $Q(x_t)$ and positive constants $\tilde \gamma_1$ and $\tilde \gamma_2$, such that $\tilde{\gamma_1} I \leq Q(x_t) \leq \tilde{\gamma_2} I$ and,
\begin{eqnarray}
E_{ \xi_t}\left[{\cal A}(x_t, \xi_t)Q(x_t){\cal A}'(x_t, \xi_t)\right] < Q(x_{t+1}).\label{Q_corollary}
\end{eqnarray}
\end{corollary}
\begin{IEEEproof}
The proof follows from Theorem \ref{lya_theorem} and by constructing $Q(x_t) = P(x_t)^{-1}$. \end{IEEEproof}
\begin{remark} \label{remark_Lyapunov}We will refer to matrix $Q(x_t)$, satisfying the conditions (\ref{Q_corollary}) of Corollary \ref{corollary_MSS} as the matrix Lyapunov function.
\end{remark}
Our goal is to derive a necessary condition for the MSE stability of the linearized error dynamics; thereby, providing a necessary condition for  MSE stability of the true error dynamics.
\begin{lemma}\label{optimal_gain_thm}
The necessary condition for exponential mean square stability of the linearized error dynamics (\ref{linearized_error_sys}) is given by
\begin{eqnarray}
(1-p)^M\left(\det(A(x_t))\right)^2\frac{\det(Q_0(x_t))}{\det(Q_0(x_{t+1}))}<1,\label{111}
\end{eqnarray}
for Lebesgue almost all $x_t\in X$. In (\ref{111}) $Q_0(x_t)$ is a solution of the following Riccati equation,
\begin{align}
Q_0(x_{t+1}) &= A(x_t)Q_0(x_t)A'(x_t) + R(x_t)\nonumber\\
& - \big[(A(x_t)Q_0(x_t)C'(x_t)\big]  \Big[\left(I_M+C(x_t)Q_0(x_t)C'(x_t)\right)^{-1}\Big]\big[C(x_t)Q_0(x_t)A'(x_t)\big],
\label{optimal_riccati}
\end{align}
where $R(x_t) \geq 0$ is some symmetric positive semi-definite matrix. Furthermore, $Q_0(x_t)$ is uniformly bounded above and below with $A(x_t) := \frac{\partial f}{\partial x}(x_t)$, $C(x_t) := \frac{\partial h}{\partial x}(x_t)$, $x_{t+1} = f(x_t)$, $I_M$ is $M\times M$ identity matrix, and $(1-p)$ is the probability of erasure.
\end{lemma}
\begin{IEEEproof}
Using the result of Corollary \ref{corollary_MSS}, the necessary condition for MSE stability of (\ref{linearized_error_sys}) can be expressed in terms of the existence of  $\tilde  \gamma_1 I\leq Q(x_t)\leq \tilde \gamma_2 I$, such that $\tilde \gamma_1, \tilde \gamma_2 >0$ and,
\begin{align}
E_{\xi_t}\left[{\cal A}(x_t,\xi_t)Q(x_t){\cal A}'(x_t,\xi_t)\right] < Q(x_{t+1}),\label{QQ}
\end{align}
where $\mathcal{A}(x_t,\xi_t) = A(x_t) - \xi_t\tilde{K}(x_t)C(x_t)$ and $\tilde{K}(x_t) := \frac{\partial K}{\partial y}(h(x_t))$. Minimizing trace of the left-hand side of (\ref{QQ}) with respect to $\tilde K(x_t)$, we obtain $\tilde K^*(x_t) = A(x_t)Q(x_t)C'(x_t)\left(C(x_t)Q(x_t)C'(x_t)\right)^{-1}$ and $Q(x_t)$ to satisfy
\begin{align}
Q(x_{t+1}) &> A(x_t)Q(x_t)A'(x_t) \nonumber\\
&\quad-pA(x_t)Q(x_t)C'(x_t)\left(C(x_t)Q(x_t)C'(x_t)\right)^{-1}C(x_t)Q(x_t)A'(x_t).
\label{eqnn1}
\end{align}
It is important to notice that the inequality (\ref{eqnn1}) is independent of any positive scaling i.e., if $Q(x_t)$ satisfies the above inequality then $c Q(x_t)$ also satisfies the above inequality for any positive constant $c$. Since $Q(x_t)$ is a matrix Lyapunov function and hence lower bounded, it follows from Remark \ref{remark_unstable}, that there exists a positive constant $\Delta>0$ such that $C(x_t)Q(x_t)C'(x_t)\frac{(1-p)}{p} \geq \Delta I_M$. Hence (\ref{eqnn1}) implies following inequality to be true
\begin{align}
Q(x_{t+1}) &> A(x_t)Q(x_t)A'(x_t)\nonumber\\
&\quad - A(x_t)Q(x_t)C'(x_t)\left(\Delta I_M+C(x_t)Q(x_t)C'(x_t)\right)^{-1}C(x_t)Q(x_t)A'(x_t).
\label{eqnn2}
\end{align}
Now define $Q_0(x_t):=\frac{1}{\Delta}Q(x_t)$, then using the fact that (\ref{eqnn2}) is independent of positive scaling, we obtain following inequality for $Q_0(x_t)$
\begin{align}
Q_0(x_{t+1}) &> A(x_t)Q_0(x_t)A'(x_t)\nonumber\\
&\quad -A(x_t)Q_0(x_t)C'(x_t)\left(I_M+C(x_t)Q_0(x_t)C'(x_t)\right)^{-1}C(x_t)Q_0(x_t)A'(x_t).
\label{eqnn3}
\end{align}
Inequality (\ref{eqnn3}) implies there exists $R(x_t)\geq 0$, such that the following equality is true.
\begin{align}
Q_0(x_{t+1})&= A(x_t)Q_0(x_t)A'(x_t) + R(x_t)\nonumber\\
&\quad-A(x_t)Q_0(x_t)C'(x_t)\left(I_M + C(x_t)Q_0(x_t)C'(x_t)\right)^{-1}C(x_t)Q_0(x_t)A'(x_t).
\label{eqnn4}
\end{align}
For any fixed trajectory $\{x_t\}$ generated by the system, $x_{t+1}=f(x_t)$, the above equality resembles the Riccati equation obtained for the minimum covariance estimator design problem for the linear time varying system, where the matrices $Q_0(x_t)$ and $R(x_t)$ can be identified with the error and input noise covariance matrices, respectively \cite{kwakernaak} with output noise variance matrix equal to identity matrix. The difference between the regular Riccati equation obtained from the minimum variance estimator problem for the linear time varying system and Eq. (\ref{eqnn4}) is that, the various matrices appearing in (\ref{eqnn4}) are parameterized by $x_t$ instead of time.  Furthermore $Q_0(x_t)$ as the solution of Riccati-like equation (\ref{eqnn4}) is both bounded above and below and is proved as follows. The system matrices $A(x_t)$ and $C(x_t)$ satisfy  Assumption \ref{assume} along any given trajectory. Hence, the linearized system, $\eta_{t+1}=A(x_t)\eta_t, \zeta_{t}=C(x_t)\eta_t$, along any fixed trajectory is uniformly completely reconstructible as defined in \cite{kwakernaak} (Definition 6.6).
It then follows from \cite{Jazwinski} (Lemmas 7.1 and 7.2) that the covariance matrix $Q_0(x_t)$ is uniformly bounded above and below for all $x\in X$. The matrix $Q_0(x_t)$ satisfies (\ref{eqnn1}) follows from the definition of $Q_0(x_t)$ ( i.e., $Q_0(x_t):=\frac{1}{\Delta}Q(x_t)$) and the fact that (\ref{eqnn1}) is independent of positive scaling. We obtain,
\begin{align}
Q_0(x_{t+1}) &> A(x_t)Q_0(x_t)A'(x_t) \nonumber\\
&\quad-pA(x_t)Q_0(x_t)C'(x_t)\left(C(x_t)Q_0(x_t)C'(x_t)\right)^{-1}C(x_t)Q_0(x_t)A'(x_t).
\label{riccati_barQ}
\end{align}
This proves that $Q_0(x_t)$ obtained as a solution of Riccati-like equation is a valid matrix Lyapunov function. To derive the required necessary condition (\ref{111}), we take determinants on both sides of (\ref{riccati_barQ}) to obtain
\begin{align}
\label{nec_eqn_int}
1 &> \det\left(I_N - pC'(x_t)\left(C(x_t)Q_0(x_t)C'(x_t)\right)^{-1}C(x_t)Q_0(x_t)\right)\left(\det(A(x_t))\right)^2\frac{\det(Q_0(x_t))}{\det(Q_0(x_{t+1}))}.
\end{align}
By Sylvester's determinant Theorem (i.e., $\det(I_N + GJ) = \det(I_M+JG)$, $G\in \mathbb{R}^{N\times M}, J\in \mathbb{R}^{M\times N}$), we obtain
\begin{align}
\label{det_thm}
(1-p)^M = \det\Big(I_N - pC'(x_t)\left(C(x_t)Q_0(x_t)C'(x_t)\right)^{-1}C(x_t)Q(x_t)\Big).
\end{align}
We obtain the required inequality (\ref{111}) by combining Eqs. (\ref{nec_eqn_int}) and (\ref{det_thm}).%
\end{IEEEproof}
The results of Lemma \ref{optimal_gain_thm} will now be used to prove the main results of the paper under various assumptions on the system dynamics.
\begin{theorem}[Linear Systems]\label{linear_thm} Let $f(x)=Ax$ with $x\in \mathbb{R}^N$ and $h(x)=Cx\in \mathbb{R}^M$. Assume that all eigenvalues $\lambda_k$ for $k=1,\ldots,N$ of $A$ have absolute value greater than one. The necessary condition for the observer error dynamics to be MSE stable is given by
\begin{align}
(1-p)^M\left(\prod_{k=1}^N|\lambda_k|\right)^2 < 1.\label{linear_necc_M}
\end{align}
\end{theorem}
\begin{IEEEproof}
For the linear system, the solution of Riccati-like equation (\ref{optimal_riccati}) from Lemma \ref{optimal_gain_thm} leads to a constant matrix $Q_0$ independent of $x_t$. Hence the necessary condition (\ref{111}) for the stability will reduce to
\[(1-p)^M\det( A^2)<1.\]
The required necessary condition (\ref{linear_necc_M}) then follows by substituting $\det(A^2)=\left(\prod_{k=1}^N |\lambda_k|\right)^2$.
\end{IEEEproof}
\begin{remark} A careful examination of the proofs for Lemma \ref{linear_error_MSS} and \ref{optimal_gain_thm}, and Theorem \ref{lya_theorem} for the special case of linear systems with single output, reveals the necessary condition (\ref{linear_necc_M}) is also sufficient for  MSE stability  of the linear system.
\end{remark}
\begin{theorem}[Nonlinear systems on unbounded space]\label{non-compact} Consider  system (\ref{dynamical_sys}) with system mapping $f$ and output $h$ satisfying Assumption \ref{assume} and state space $X$ possibly unbounded. The necessary condition for  MSE stability of the observer error dynamics (\ref{observer_sys}) is given by
\begin{eqnarray}
(1-p)^M\left(\det(A(x_t))\right)^2\frac{\det(Q_0(x_t))}{\det(Q_0(x_{t+1}))}<1,\label{necessary_unbounded}
\end{eqnarray}
for Lebesgue almost all $x\in X$, where $A(x)=\frac{\partial f}{\partial x}(x)$ and $Q_0(x)$ satisfy the Riccati-like Eq. (\ref{optimal_riccati}).
\end{theorem}
\begin{IEEEproof}
The proof follows by combining results from Lemmas \ref{linear_error_MSS} and \ref{optimal_gain_thm}, and Theorem \ref{lya_theorem}.
\end{IEEEproof}
In Theorem \ref{compact}, we show, for a nonlinear system evolving on a compact state space, the term $\left(\det(A(x_t))\right)^2\frac{\det(Q(x_t))}{\det(Q(x_{t+1}))}$ from (\ref{necessary_unbounded}) relates to the sum of postive Lyapunov exponents of the system.
For Theorem \ref{compact} we provide the following definitions \cite{Mane}.
\begin{definition}[Physical measure]\label{invariant_measure} Let ${\cal M}(X)$ be the space of probability measures on $X$. A measure $\mu\in {\cal M}(X)$ is said to be invariant for $x_{t+1}=f(x_t)$ if
$\mu(f^{-1}(B))=\mu(B)$
for all sets $B \in \mathcal{B}(X)$ (Borel $\sigma$-algebra generated by $X$). An invariant probability measure, $\mu$, is said to be ergodic if any continuous bounded function $\varphi$ that is invariant under $f$, i.e.,  $\varphi(f(x))=\varphi(x)$, is $\mu$ almost everywhere constant. Ergodic invariant measure, $\mu$, is said to be physical if $\lim_{n\to \infty}\frac{1}{n}\sum_{k=0}^n \varphi(f^k(x))=\int_X \varphi(x)d\mu(x)$ for positive Lebesgue measure of the initial condition $x\in X$ and all continuous function $\varphi: X\to \mathbb{R}$.
\end{definition}
\begin{definition}[Lyapunov exponents]\label{Lyapunov_exponents} For a deterministic system $x_{t+1}=f(x_t)$, let
\begin{align}
\Lambda(x_0)=\lim_{t\to \infty}\left(D_x^tf(x_0)'D_x^tf(x_0)\right)^{\frac{1}{2t}},\label{le_average}
\end{align}
where $D_x f(x)=\frac{\partial f}{\partial x}(x)$ and $D_x^t f(x_0):=D_x f(x_t)\cdots D_xf(x_0)$. Let $\lambda^i_{exp}$ for $i=1,\ldots,N$ be the eigenvalues of $\Lambda(x_0)$, such that $\lambda^1_{exp} \geq \lambda^2_{exp} \geq \cdots \geq \lambda^N_{exp}$. Then, the Lyapunov exponents  $\Lambda^i_{exp}$ are defined as $\Lambda^i_{exp}=\log \lambda^i_{exp}$ for $i=1,\ldots,N$. Furthermore, if $\det\left(\Lambda(x_0)\neq 0\right)$, then
\begin{align}
\label{det_sum_exp}
\lim_{t\to \infty}\frac{1}{t}\log\left|\det\left(D_x^tf(x_0)\right)\right| =\log \prod_{k=1}^N \lambda^k_{exp}(x).
\end{align}
\end{definition}
\begin{remark}\label{remark_ergodic}
The technical conditions for the existence of limits in (\ref{le_average}) and (\ref{det_sum_exp}) are provided by  the Multiplicative Ergodic Theorem \cite{Ruelle_ergodic} (Theorem 1.6), \cite{walter_ergodic_theory} (Theorem 10.4),  \cite{Ruelle85} (Section D). The limits in (\ref{le_average}) and (\ref{det_sum_exp}) are  known to be independent of the initial condition and are unique under the assumption of unique ergodic invariant measure for system dynamics. For a compact state space, the existence of an invariant measure is always guaranteed \cite{walter_ergodic_theory} (Corollary 6.9.1). Furthermore, every invariant measure admits ergodic decomposition \cite{walter_ergodic_theory} (Remarks pp. 153), \cite{Mane} (Theorem 6.4).
We now make Assumption \ref{unique} on the system dynamics.
\end{remark}
\begin{assumption}\label{unique}
We assume  the nonlinear system, $x_{t+1}=f(x_t)$, has a unique physical measure with all Lyapunov exponents positive.
\end{assumption}
The assumption of a unique physical measure is not restrictive and it allows us to prove the main result in Theorem \ref{compact}, that is independent of initial conditions. With ergodic invariant measures that are guaranteed to exist (Remark \ref{remark_ergodic}), the main result in Theorem \ref{compact} will be a function of a particular ergodic measure under consideration. The assumption of all Lyapunov exponent being positive is analogous to the assumption made in the LTI case that all eigenvalues are positive. We verify through simulation results in section \ref{sim} that the result of Theorem \ref{compact} also applies to the case where one of the Lyapunov exponent is negative.

\begin{theorem}[Nonlinear systems on compact space] \label{compact}
Consider the system (\ref{dynamical_sys}) with system mapping $f$ and output $h$ satisfying Assumptions \ref{assume} and \ref{unique} and state space $X$ compact. The necessary condition for  MSE stability of the observer error dynamics (\ref{observer_sys}) is given by
\begin{align}
(1-p)^M\left(\prod_{k=1}^N \lambda^k_{exp}\right)^2<1, \label{condition_M}
\end{align}
where $\lambda_{exp}^k = e^{\Lambda_{exp}^k}$, and $\Lambda_{exp}^k$ is the $k^{th}$  positive Lyapunov exponent of  $x_{t+1} = f(x_t)$.
\end{theorem}
\begin{IEEEproof} We follow the notations from Lemma \ref{optimal_gain_thm}. The necessary condition for  MSE stability (Eq. \ref{111}) is true for almost all points $x\in X$, and, hence in particular for $x_t$ evaluated along the system trajectory $x_{t+1}=f(x_t)$. Evaluating (\ref{111}) along the system trajectory and taking the product, we write the necessary condition as
\[
\left((1-p)^M\right)^n \det (Q_0(x_0)Q_0^{-1}(x_{n+1}))\prod_{t=1}^n \det (A(x_t))^2<1.\]
Taking time average for the $\log$ of the expression and in the limit as $n\to \infty$, we obtain the following necessary condition for MSE stability,
\begin{align}
\lim_{n\to \infty}&\frac{1}{n}\Bigg[\log \left((1-p)^M\right)^n + \log\left(\det (Q_0(x_0)Q_0^{-1}(x_{n+1}))\prod_{t=1}^n \det (A(x_t))^2\right)\Bigg] <0.\label{ss}
\end{align}
Using the fact that both $Q_0(x_t)$ and $Q_0^{-1}(x_t)$ are almost always uniformly bounded and using (\ref{det_sum_exp}) from Definition \ref{Lyapunov_exponents}, (\ref{ss}) gives the required necessary condition (\ref{condition_M}) for MSE stability.
\end{IEEEproof}
\begin{remark} The necessary condition for MSE stability  in Theorems \ref{linear_thm}, \ref{non-compact}, and \ref{compact} for single input case is tighter however for $1<M<N$, we expect the condition to be improved further.
The necessary condition for MSE stability from our main  results provides a critical dropout rate, i.e., the erasure probability, $q^{*}=1-p^{*}$, above which the system is guaranteed  MSE unstable. In particular, the critical dropout rate for a nonlinear system with single output, evolving  on compact space from Theorem \ref{compact} is given by $q^{*}=\left(\prod_{k=1}^N \lambda^k_{exp}\right)^{-2}$.
\end{remark}
\subsection{Entropy and limitation for observation}
Measure-theoretic entropy, $H_\mu(f)$, for the dynamical system, $x_{n+1}=f(x_n)$, is associated with a particular ergodic invariant measure, $\mu$, and is another measure of dynamical complexity. While the measure-theoretic entropy counts the number of {\it typical} trajectories for their growth rate, the positive Lyapunov exponents measure the rate of exponential divergence of nearby system trajectories. For more details on entropy refer to \cite{walter_ergodic_theory}. These two measures of dynamical complexity are related by Ruelle's inequality.
\begin{theorem}[Ruelle's Inequality] \label{theorem_ruelle_inequality}(\cite{Ruelle85} Eq. 4.4); (\cite{Ruelle_inequality} Theorem 2) Let $x_{n+1}=f(x_n)$ be the dynamical system, $f:X\to X$ be a $C^r$ map, with $r\geq 1$, of a compact metric space $X$ and $\mu$ an ergodic invariant measure. Then, 
\begin{eqnarray}
H_\mu(f)\leq\sum_{k} (\Lambda^k_{exp})^{+}, \label{ruellle_inequality}
\end{eqnarray}
where $a^{+}=\max\{0,a\}$, $H_\mu(f)$ is the measure-theoretic entropy corresponding to the ergodic invariant measure $\mu$, and $\Lambda^{k}_{exp}$ are the Lyapunov exponents of the system.
\end{theorem}
The Ruelle inequality (\ref{ruellle_inequality}) can be used to relate the limitation for observation with system entropy.
\begin{theorem}\label{entropy_theorem} Consider the system (\ref{dynamical_sys}) with system mapping $f$ and output $h$ satisfying Assumptions \ref{assume} and \ref{unique} and state space $X$ compact. The necessary condition for  MSE stability of the observer error dynamics (\ref{observer_sys}) is given by
\begin{align}
M\log (1-p)+2 H_\mu(f) < 0
\end{align}
where $\mu$ is the physical invariant measure of $f$ (Definition \ref{invariant_measure} and Assumption \ref{unique}) and $H_\mu(f)$ is the measure-theoretic entropy corresponding to measure $\mu$.
\end{theorem}
\begin{IEEEproof} The proof follows by applying the results of Theorems \ref{compact}
 and \ref{theorem_ruelle_inequality}.
 \end{IEEEproof}
\section{Simulation Results}
\label{sim}
\begin{figure}[ht!]
\begin{center}
\mbox{
\hspace{-0.1in}
\subfigure[]{\scalebox{0.34}{\includegraphics{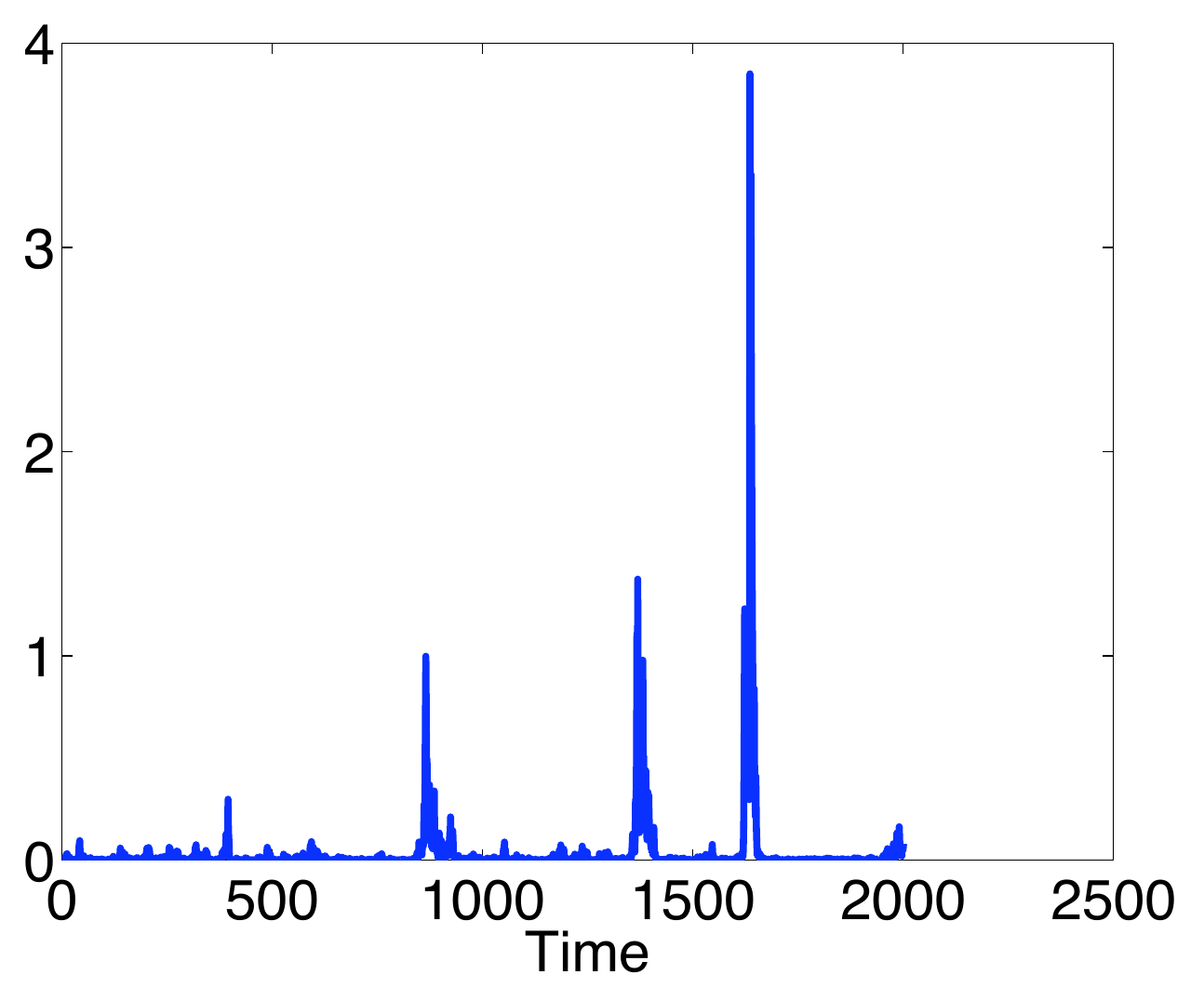}}}
\hspace{-0.14in}
\subfigure[]{\scalebox{0.25}{\includegraphics{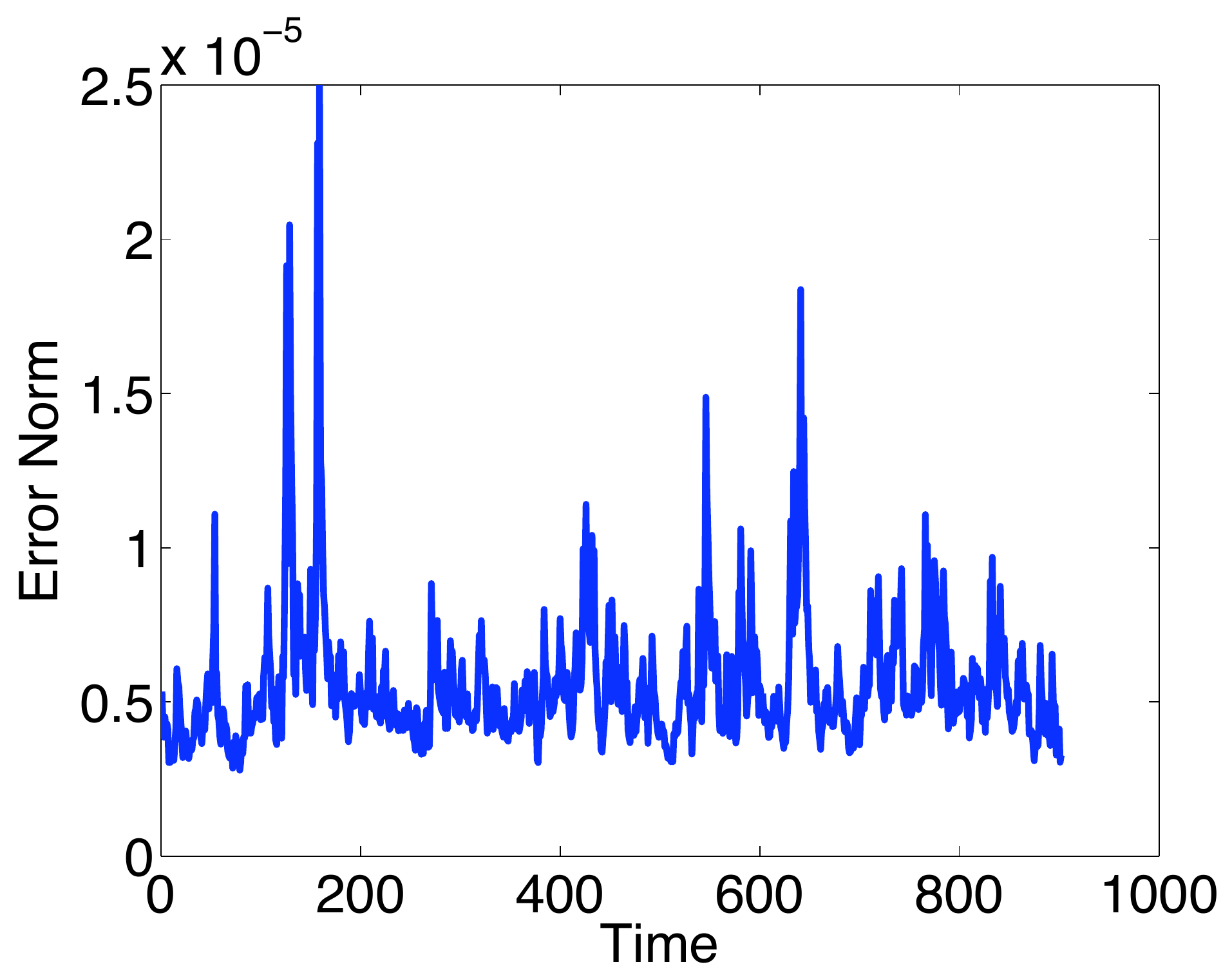}}}}
\mbox{
\hspace{-0.12in}
\subfigure[]{\scalebox{0.28}{\includegraphics{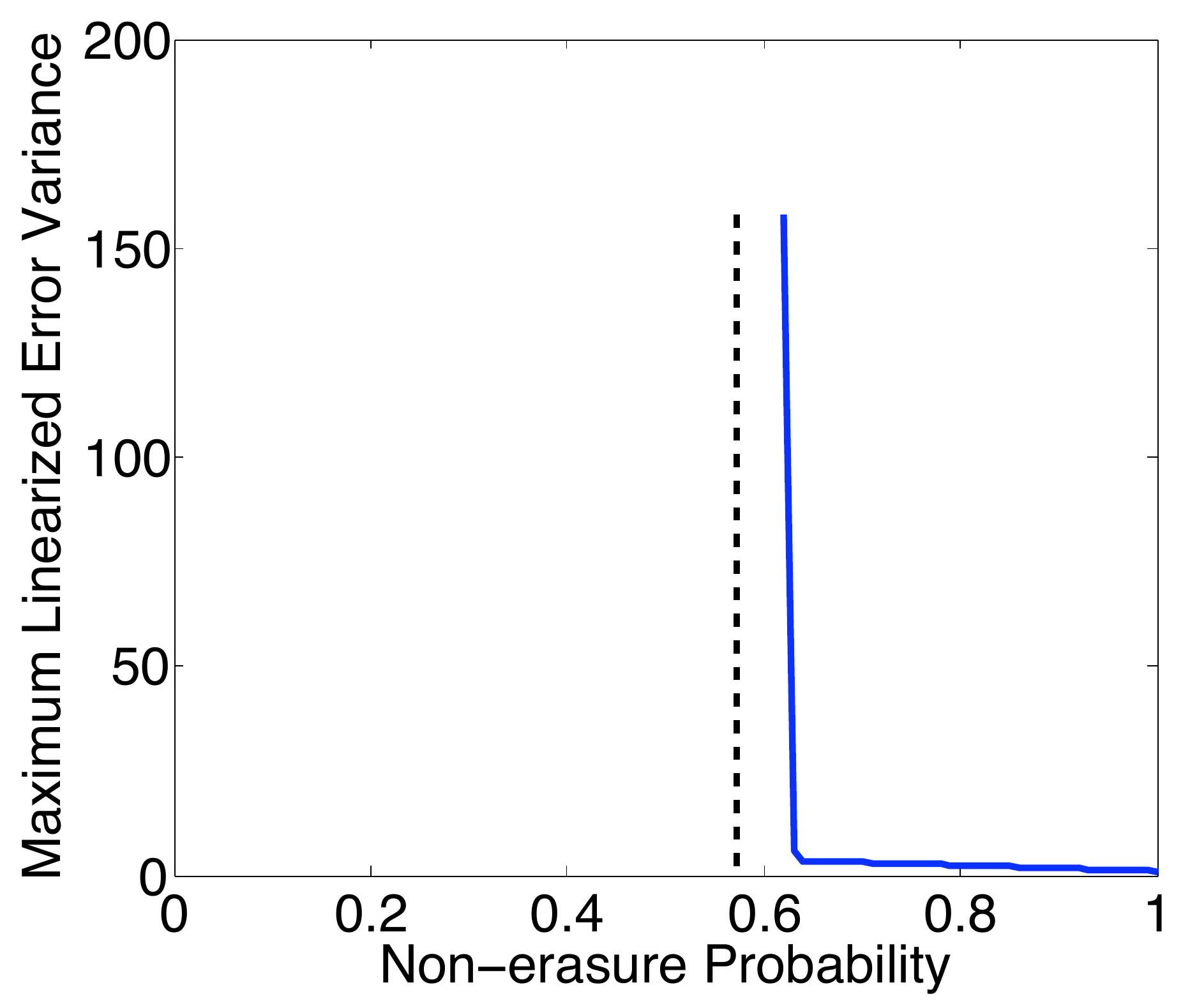}}} }
 \caption{(a) Error norm as a function of time for $p=0.55$; (b) Error norm as a function of time for $p=0.7$; (c) Maximum linearized covariance vs non-erasure probability for Henon map}
\label{nonlinear_obsvr}
\end{center}
\end{figure}

\noindent {\bf Henon map} is one of the widely studied examples of two-dimensional chaotic maps. The small random perturbation of a two-dimensional Henon map is described by following equations:
\begin{align*}
x_{1t+1} &= 1 - ax_{1t}^2 + x_{2t} + r_{1t},\\
x_{2t+1}&=b x_{1t} + r_{2t},\\
y_t &= \xi_t x_{1t},
\end{align*}
where  $a = 1.4$, $b = 0.3$ are constant parameters, and $r_{it}\in [0,1\text{E-6}],\; i\in \{1,2\}$, are uniform random variables. The small amount of external noise, $r_{it}$, is essential to see the effect of mean square instability. The system has Lyapunov exponents given by $\lambda_1 = 0.426$ and $\lambda_2 = -1.63$. Although the main results of this paper are proved under the assumption that all Lyapunov exponents are positive, the simulation results verify that the results hold true even for this example with one Lyapunov exponent negative. The critical probability $p^{*}$ is computed, based on the positive exponent and is equal to $p^{*}=1-\frac{1}{\exp^{2\lambda_1}}=0.5734$. The observer is designed such that error dynamics with no erasure is asymptotically stable. In Figs. (\ref{nonlinear_obsvr}a) and (\ref{nonlinear_obsvr}b), we plot the error norm for the observer dynamics, averaged over $50$ realizations of the erasure sequence, at probabilities below and above the critical probability $p^*$, respectively. We clearly see the average error norm for non-erasure probability, $p=0.7>p^*$, is negligible compared to fluctuations in the average error norm for $p=0.55<p^*$, which are four orders of magnitude higher than the uniform noise in the system. In Fig. (\ref{nonlinear_obsvr}c), we plot the peak error variance for linearized error dynamics vs. non-erasure probability. The dashed line indicates the critical probability, $p^* = 0.5734$. We observe the peak linearized error variance is unbounded below critical probability.
\section{Conclusions}
\label{inf}
In this work, the problem of state observation for a nonlinear system over erasure channel is studied. The main results of this paper prove that limitation arises for MSE stabilization of observer error dynamics. We show that instability of the non-equilibrium dynamics of the nonlinear system, as captured by positive Lyapunov exponents, plays an important role in obtaining the limitation result for nonlinear observation. The limitation result for LTI systems is obtained as a special case, where Lyapunov exponents emerge as the natural generalization of eigenvalues from linear systems to nonlinear systems. The proof technique presented in this paper can be easily extended to prove results for the estimation of linear time varying systems over erasure channels.
\section{Acknowledgment}The research work was supported by National Science Foundation (CMMI  0807666) and (ECCS 1002053) grant. The authors would like to thank Prof. Nicola Elia for useful discussion.

\bibliographystyle{IEEEtran}
\bibliography{ref1,ref}
\end{document}